\documentclass[12pt,reqno]{amsart}
\usepackage{latexsym,amsmath}
\usepackage[left=3cm,top=2.5cm,right=3cm,bottom=2.5cm]{geometry}
\usepackage{color}
\usepackage{amsthm}
\usepackage{amssymb}
\usepackage{epsfig}
\begin{document}

\newtheorem{theorem}{Theorem}[section]
\newtheorem{cor}[theorem]{Corollary}
\newtheorem{lemma}[theorem]{Lemma}
\newtheorem{fact}[theorem]{Fact}
\newtheorem{property}[theorem]{Property}
\newtheorem{proposition}[theorem]{Proposition}
\newtheorem{claim}[theorem]{Claim}
\newtheorem{definition}[theorem]{Definition}
\theoremstyle{definition}
\newtheorem{example}[theorem]{Example}
\newtheorem{remark}[theorem]{Remark}
\newcommand\eps{\varepsilon}
\newcommand\la {\lambda}
\newcommand{\E}{\mathbb E}
\newcommand{\Var}{{\rm Var}}
\newcommand{\Prob}{\mathbb{P}}
\newcommand{\N}{{\mathbb N}}
\newcommand{\eqn}[1]{(\ref{#1})}
\newcommand{\gP}{\mathcal{P}}
\newcommand{\gQ}{\mathcal{Q}}
\newcommand{\dist}{d}

\def\G{ \mathcal{G}}
\def\pr{\Prob}

\title{Geodesics and almost geodesic cycles in random regular graphs}
\author{Itai Benjamini}
\address{Department of Mathematics\\
The Weizmann Institute of Science \\
Rehovot 76100 \\
Israel} \email{\tt itai.benjamini@weizmann.ac.il}
\author{Carlos Hoppen}
\address{Department of Combinatorics and Optimization \\
University of Waterloo \\
Waterloo-ON \\
Canada N2L 3G1} \email{\tt choppen@math.uwaterloo.ca}
\author{Eran Ofek}
\address{Department of Computer Science and Applied Mathematics \\
The Weizmann Institute of Science \\
Rehovot 76100 \\
Israel} \email{\tt eran.ofek@weizmann.ac.il}
\author{Pawe{\l} Pra{\l}at}
\address{Department of Mathematics and Statistics \\
Dalhousie University \\
Halifax-NS \\
Canada B3H 3J5} \email{\tt pralat@mathstat.dal.ca}
\author{Nick Wormald} \thanks{The fifth author acknowledges the support of
the  Canadian Research Chairs Program and NSERC}
\address{Department of Combinatorics and Optimization \\
University of Waterloo \\
Waterloo-ON \\
Canada N2L 3G1}
\email{\tt nwormald@math.uwaterloo.ca}
\date{}

\begin{abstract}
A geodesic in a graph $G$ is a shortest path between two vertices of
$G$. For a specific function $e(n)$ of $n$, we define an almost
geodesic cycle $C$ in $G$ to be a cycle in which for every two
vertices $u$ and $v$ in $C$, the distance $d_G(u,v)$ is at least
$d_C(u,v)-e(n)$. Let $\omega(n)$ be any function tending to infinity
with $n$. We consider a random $d$-regular graph on $n$ vertices. We
show that almost all pairs of vertices belong to an almost geodesic
cycle $C$ with $e(n)= \log_{d-1} \log_{d-1} n +\omega(n)$ and
$|C|=2\log_{d-1}n+O(\omega(n))$. Along the way, we obtain results on
near-geodesic paths. We also give the limiting distribution of the
number of geodesics between two random vertices in this random
graph.
\end{abstract}

\maketitle

\section{Introduction}

A {\em geodesic} in a graph $G$ is a shortest path between two
vertices of $G$. Let $\omega(n)$ be any function tending to infinity
with $n$,   and put  $e(n)= \log_{d-1} \log_{d-1} n +\omega(n)$. We
define an {\em almost geodesic cycle} $C$ in $G$ to be a cycle in
which for every two vertices $u$ and $v$ in $C$, the distance
$d_G(u,v)$ is at least $d_C(u,v)-e(n)$.  We investigate the
existence of almost geodesic cycles through random pairs of vertices
in a random graph, and related questions on geodesics and paths that
are  nearly geodesic, in a sense to be made precise. Our results
refer to the probability space of random $d$-regular graphs with
uniform probability distribution. This space is denoted
$\mathcal{G}_{n,d}$, and asymptotics (such as ``asymptotically
almost surely'', which we abbreviate to a.a.s.) are for $n\to\infty$
with $d\ge 3$ fixed, and $n$ even if $d$ is odd.

Some related previous research focussed on finding
(edge/internally)-disjoint paths with many sources and targets.
Frieze and Zhao~\cite{FZ} showed that for sufficiently large $d$
there exist fixed positive constants $\alpha$ and $\beta$ such that
a graph $G$ taken from $\mathcal{G}_{n,d}$ a.a.s.\ has the following
property: for any choice of $k$ pairs $\{(a_i,b_i) ~|~
i=1,\ldots,k\}$, satisfying
\begin{itemize}
\item [(i)] $k \leq \lceil \alpha d n / \log_d n \rceil$, and
\item [(ii)] for each vertex $v$: $|i: a_i =v| + |i: b_i =v| \le \beta d$,
\end{itemize}
there exist edge-disjoint paths in $G$ connecting $a_i$ to $b_i$ for
all $i=1,2,\ldots, k$.   This result is optimal up to constant
factors. The paths returned by their algorithm are of length of at
least $10 \log_d n$.

Our focus is different as it comes from different motivation:
studying almost geodesic cycles in $\mathcal{G}_{n,d}$. Our result
on internally disjoint paths refers to one pair of vertices fixed
before the graph is chosen. This is a much weaker model than the
model of~\cite{FZ}, that dealt with $\Theta(n/ \log n)$ pairs given
by an adversary after the graph is chosen. However, we show the
existence of disjoint paths that approximate the optimal path (whose
length is a.a.s.\ in $[\log_{d-1} n -\omega (n), \log_{d-1} n +
\omega(n)]$) by an additive factor of $\log_{d-1} \log_{d-1} n$,
whereas the result of~\cite{FZ} give at best a constant
multiplicative factor. Additionally, our result holds for all $d \ge
3$, and that we find the maximum possible number of internally
disjoint paths, $d$, that there can possibly be between two
vertices.

\begin{theorem}\label{main}
Take any integer $d \ge 3$ and any function $\omega(n)$ with
$\omega(n) \rightarrow \infty$. Let $G \in\mathcal{G}_{n,d}$ and
choose vertices $u$ and $v$ in $V(G)$ independently with uniform
probability. Then a.a.s.\ the following hold:
\begin{itemize}
\item[(i)]  $|\dist(u,v) - \log_{d-1}n| <\omega(n)$,
\item[(ii)]  there are $d$ paths connecting
$u$ and $v$ such that the subgraph induced by each pair of these
paths is an almost geodesic cycle.
\end{itemize}
\end{theorem}
Note that the $d$ paths in (ii) theorem  are internally disjoint
because each pair of them induces a cycle.

In a slightly different direction, we also investigate the
distribution of the number of geodesics joining two vertices (see
Theorem~\ref{thm3.5}).

We may obtain the lower bound in part~(i) of the theorem from an
elementary observation. Note that, given $G \in \mathcal{G}_{n,d}$,
the number of vertices at distance at most $i$ from a vertex $u$ is
bounded above by
$$1+d+d(d-1)+\ldots+d(d-1)^{i-1}=O\left( (d-1)^i \right).$$
So, there are $O\left( n (d-1)^{-\omega(n)} \right)$ vertices at
distance $i=\log_{d-1}n - \omega(n)$ from any given vertex of $G$,
where $\omega(n) \rightarrow \infty$. As a consequence, if two
vertices of $G$ are chosen independently with uniform probability,
then the probability that the second vertex is at distance at most
$i=\log_{d-1}n - \omega(n)$ from the first is at most
$$\frac{1}{n} O\left( n (d-1)^{-\omega(n)} \right) =  O\left((d-1)^{-\omega(n)} \right),$$
and therefore, a.a.s.\
\begin{equation} \label{l:first_time}
d(u,v) \geq \log_{d-1}n - \omega(n)
\end{equation}
if $u,v$ are vertices chosen independently with uniform
probability in $G \in \mathcal{G}_{n,d}$ and $\omega(n)$ is a
function satisfying $\omega(n) \rightarrow \infty$. The fact that
a.a.s.\ $d(u,v) \leq \log_{d-1}n + \omega(n)$ will follow from our
study of the distribution of the number of geodesics in $G \in
\mathcal{G}_{n,d}$.

The rest of the proof requires more sophisticated arguments. Instead
of working directly in the uniform probability space of random
regular graphs on $n$ vertices $\mathcal{G}_{n,d}$, we use the
\textit{pairing model} of random regular graphs, first introduced by
Bollob\'{a}s~\cite{bollobas2}, which is described next. Suppose that
$dn$ is even, as in the case of random regular graphs, and consider
$dn$ points partitioned into $n$ labelled cells $v_1,\ldots,v_n$ of
$d$ points each. A \textit{pairing} of these points is a perfect
matching of them into $dn/2$ pairs. Given a pairing $P$, we may
construct a multigraph $G(P)$, with loops allowed, as follows: the
vertices are the cells $v_1$,\ldots, $v_n$, and a pair $\{x,y\}$ in
$P$ corresponds to an edge $v_iv_j$ in $G(P)$ if $x$ and $y$ are
contained in the cells $v_i$ and $v_j$, respectively. It is an easy
fact that the probability of a random pairing corresponding to a
given simple graph $G$ is independent of the graph, hence the
restriction of the probability space of random pairings to simple
graphs is precisely $\mathcal{G}_{n,d}$. Moreover, it is well known
that a random pairing generates a simple graph with probability
asymptotic to a constant depending on $d$, so that any event holding
a.a.s.\ over a probability space of random pairings also holds
a.a.s.\ over the corresponding space $\mathcal{G}_{n,d}$. For this
reason, asymptotic results over random pairings suffice for our
purposes. The advantage of using this model is that the pairs may be
chosen sequentially so that the next pair is chosen uniformly at
random over the remaining (unchosen) points. For more information on
this model, see~\cite{nick1}.

The numbers of geodesics are investigated in Section~\ref{section3}.
Theorem~\ref{main} is proved in Section~\ref{almost}. Some final
remarks are in Section~\ref{final}.

\section{Distribution of the number of geodesics}\label{section3}

The first portion of our argument is a simplified version of part of
the argument of Bollob{\'a}s and Fernandez de la Vega~\cite{BF}. We
consider the process in which the neighbourhoods of $u$ and $v$ are
exposed step by step. First, the neighbours of $u$ and $v$ are
revealed, then the vertices at distance two, and so on.  This
sequential exposure of the random regular graph is analysed using
the random pairing model mentioned in the Introduction.

Let $N_i(u)$ denote the set of vertices at distance at most $i$ from
$u$. Note that, in the early stages of this process, the graphs
grown from $u$ and $v$ tend to be trees, hence the number $n_i$ of
elements in $N_i(u)$ is approximately
$$n_{i-1}+(d-1)(n_{i-1}-n_{i-2})\,.$$

Let $f_i$ denote the number of vertices in a balanced $d$-regular
tree, that is,
$$f_i = 1 + d \sum_{j=0}^{i-1}(d-1)^j = 1 +
\frac{d((d-1)^{i}-1)}{d-2}\,,$$ and let
$$ i_0 = \left\lfloor \frac 12 \log_{d-1}n \right\rfloor \,.$$

\begin{lemma}\label{lem:number_of_nodes}
Let $\omega(n)$ be any function of $n$ such that $\omega(n)
\rightarrow \infty$. For $i \leq i_0 - \omega(n)$ a.a.s.\ the
cardinality $n_i$ of $N_i(u)$ equals $f_i$. Moreover, for $i \leq
i_0 + \omega(n)$ a.a.s.\
$$ n_i = f_i - O\Big(\omega(n)(d-1)^{3(i-i_0)+\omega(n)}\Big) \,.$$
\end{lemma}

\begin{proof}
First note that it is sufficient to consider the case when
$\omega(n)=o(\log n)$.

Since $f_i$ denotes the number of vertices in a balanced tree where
every non-leaf vertex has degree $d$, the first assertion follows if
we show that a.a.s.\ the set of vertices at distance at most $i \leq
i_0 - \omega(n)$ of a vertex $u$ induces a tree. In other words, if
we expose, step by step, the vertices at distance $1,2,\dots, i$
from $u$, we have to avoid, at step $j$, edges that induce cycles.
So, we wish not to find edges between any two vertices at distance
$j$ from $u$ or edges that join any two vertices at distance $j$ to
a same vertex at distance $j+1$ from $u$. We shall refer to edges of
this form as `bad edges'. Note that the expected number of `bad
edges' at step $i+1$ is equal to
$O(n_i^2/n)=O(f_i^2/n)=O((d-1)^{2i}/n)$.

Consider $i_1 = \lfloor \frac 12 \log_{d-1}n - \omega(n) \rfloor$.
The expected number of `bad edges' found up to step $i_1$ is equal
to
$$\sum_{j=0}^{i_1-1} O\big((d-1)^{2j} / n\big)  = O\big((d-1)^{2i_1}
/ n\big) = O\big((d-1)^{-2\omega(n)}\big) = o(1)\,.$$ Thus, by
Markov's inequality, a.a.s.\ there are no `bad edges' until step
$i_1$, hence a.a.s.\ $N_{i_1}(u)$ is a tree and $n_i=f_i$ for $i
\leq i_1$.

In order to prove the second assertion, notice that the expected
number of `bad edges' added between step $i_1+1$ and step $i$, $i
\le \lfloor i_0 + \omega(n) \rfloor = \lfloor \frac 12 \log_{d-1} n
+ \omega(n) \rfloor$ is equal to
$$\sum_{j=i_1}^{i-1} O\big((d-1)^{2j} / n\big) = O\big((d-1)^{2i} /
n\big) = O\big((d-1)^{2(i-i_0)}\big) \,.$$ Thus, again by Markov's
inequality, a.a.s.\ the total number of `bad edges' at time $i$ is
at most $O\big(\omega(n)(d-1)^{2(i-i_0)}\big)$. Notice that one
`bad edge' added in this time interval can destroy a tree branch
of size $O\big((d-1)^{i-i_0+\omega(n)}\big)$. This occurs because
the `bad edge' creates a cycle instead of exposing a new vertex
$v$. The branch of descendants of $v$, which would appear had $v$
been exposed and had the process continued as a balanced
$d$-regular tree, is therefore destroyed and has size at most
$1+(d-1)+\cdots+(d-1)^{i-i_0+\omega(n)}=
O\big((d-1)^{i-i_0+\omega(n)}\big)$.

Thus, we have a.a.s.\
$$n_i = f_i - O\big(\omega(n)(d-1)^{2(i-i_0)}\big) \cdot
O\big((d-1)^{i-i_0+\omega(n)}\big).$$ This completes the proof of
the lemma.
\end{proof}

Immediately from this lemma, we have
\begin{cor} \label{wasntused}
For $i = i_0 + o(\log n)$ a.a.s.\
$$ n_i = f_i - n^{o(1)} = n^{1/2 + o(1)}\,.$$
\end{cor}

In the remainder of this notes, let $u,v$ be vertices chosen
independently with uniform probability in a graph $G \in
\mathcal{G}_{n,d}$ and consider the process of exposing the
neighbourhoods of $u$ and $v$ introduced in Lemma
\ref{lem:number_of_nodes}. We say that $N_i(u)$ and $N_i(v)$
\textit{join at time $i$} if $N_{i-1}(u) \cap N_{i-1}(v) =
\emptyset$ and $N_i(u) \cap N_i(v) \neq \emptyset$.

Also, whenever a result that holds for any $\omega(n)$ satisfying
$\lim_{n \rightarrow \infty} \omega(n) = \infty$ is proven, we shall
assume without loss of generality that $\omega(n) = o (\log n)$.

\begin{lemma}\label{no join}
Let $k$ be a fixed integer and define $\gamma(n,d) = \frac{1}{2}
\log_{d-1}n - i_0$, the fractional part of $\frac{1}{2} \log_{d-1}
n$. Then
$$ \Prob (N_{i_0 + k}(u) \cap N_{i_0 + k}(v) = \emptyset) \sim
\exp \left(-\frac{d(d-1)^{2k-2 \gamma(n,d)}}{d-2}\right).$$
\end{lemma}

\begin{proof}
Denote $S_i$ the event that the neighbourhoods of $u$ and $v$ are
separate at time $i$, that is, $N_j(u)$ and $N_j(v)$ did not join up
to time $i$. We claim that
$$ \Prob (S_{i_0+k} ~ |~ S_{i_0+k-1} ) \sim \exp \left(- d^2(d-1)^{2k-2
\gamma(n,d)-2}\right).$$ This implies our result for the following
reasons. If $M$ is a positive integer, $-M<k$,
\begin{eqnarray*}
\Prob (S_{i_0+k}) & = & \Prob
(S_{i_0-M}) \\
&& \times \prod_{l=-M+1}^k \Prob (S_{i_0+l} ~ | ~ S_{i_0+l-1}).
\end{eqnarray*}
Furthermore, equation (\ref{l:first_time}) establishes that
a.a.s.\ $\dist(u,v)> 2i_0 - \omega(n)$ for any function
$\omega(n)$ with $\lim_{n \rightarrow \infty} \omega(n) = \infty$.
In particular, given $\epsilon>0$, we can choose
$M=M_{\epsilon}>0$ sufficiently large so that $\displaystyle{\Prob
(S_{i_0-M}) > 1 - \epsilon}$. Given such an $M$, we use the
previous equation to derive
\begin{equation*}
\begin{split}
\Prob (S_{i_0+k})&> (1 - \epsilon) \prod_{l=-M+1}^{k} \exp
\left(-d^2(d-1)^{2l-2
\gamma(n,d)-2}\right) \left( 1 - o(1) \right)\\
& \sim (1 - \epsilon) \exp \left( \sum_{l=-M+1}^{k}
-\frac{d^2(d-1)^{2l}}{(d-1)^{2+2\gamma(n,d)}} \right)\\
& = (1- \epsilon) \exp \left( - \frac{d\big( (d-1)^{2k+2M} -
1\big)}{(d-1)^{2M+2\gamma(n,d)}(d-2)} \right).
\end{split}
\end{equation*}
The same calculations also lead us to
\begin{equation*}
\begin{split}
\Prob (S_{i_0+k}) &< \prod_{l=-M+1}^{k} \exp
\left(-\frac{d^2(d-1)^{2l}}{(d-1)^{2+2\gamma(n,d)}}\right) \left(
1 -
o(1) \right)\\
& \sim \exp \left( -\frac{d\big( (d-1)^{2k+2M} -
1\big)}{(d-1)^{2M+2\gamma(n,d)}(d-2)} \right).
\end{split}
\end{equation*}
Putting the last two equations together and letting $\epsilon \to
0$, during which we may assume $M_{\epsilon} \rightarrow \infty$,
we have
$$ \Prob (S_{i_0+k}) \sim \exp
\left(-\frac{d(d-1)^{2k- 2\gamma(n,d)}}{d-2}\right),$$ as
required.

We now focus on proving the claim. First we would like to find the
expected number of joins at time $i=i_0+k$ given that  $N_{i-1}(u)
\cap N_{i-1}(v) = \emptyset$. Let $U_{i-1} = N_{i-1}(u) \setminus
N_{i-2}(u)$ and $V_{i-1} = N_{i-1}(v) \setminus N_{i-2}(v)$. These
are the sets of vertices at distance $i-1$ from $u$ and $v$,
respectively. We have to consider two types of join at time $i$.
The first type (see Figure~\ref{fig:first}) consists of edges that
join one vertex in $U_{i-1}$ to a vertex in $V_{i-1}$ and
therefore create $uv$-paths of odd length. The second type (see
Figure~\ref{fig:second}) contains joins such that a vertex in
$V(G) \setminus (N_{i-1}(u) \cup N_{i-1}(v))$ has neighbours in
each of $U_{i-1}$ and $V_{i-1}$. This generates a path of even
length between $u$ and $v$.

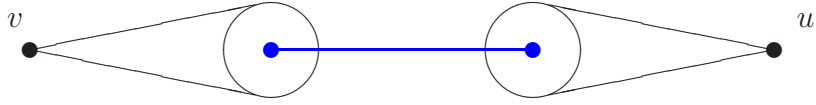
\begin{figure} \begin{center} \setlength{\unitlength}{0.3cm}
\begin{picture}(35,4) \put(1,2) {\circle*{0.7}} \put(0,3) {$v$} \put(1,2)
{\line(5,1) {10}} \put(1,2) {\line(5,-1) {10}} \put(11.7,2)
{\circle{4}} \put(34,2) {\circle*{0.7}} \put(35,3) {$u$}
\put(34,2) {\line(-5,1) {10}} \put(34,2) {\line(-5,-1) {10}}
\put(23.3,2) {\circle{4}} \color{blue} \thicklines \put(11.7,2)
{\circle*{0.7}} \put(23.3,2) {\circle*{0.7}} \put(11.7,2)
{\line(1,0) {11.6}} \end{picture} \end{center}
    \caption{First case -- odd length.} \label{fig:first} \end{figure}
\begin{figure} \begin{center} \setlength{\unitlength}{0.3cm}
\begin{picture}(35,4) \put(1,2) {\circle*{0.7}} \put(0,3) {$v$} \put(1,2)
{\line(5,1) {10}} \put(1,2) {\line(5,-1) {10}} \put(11.7,2)
{\circle{4}} \put(34,2) {\circle*{0.7}} \put(35,3) {$u$}
\put(34,2) {\line(-5,1) {10}} \put(34,2) {\line(-5,-1) {10}}
\put(23.3,2) {\circle{4}} \color{blue} \thicklines \put(11.7,2)
{\circle*{0.7}} \put(23.3,2) {\circle*{0.7}} \put(11.7,2)
{\line(1,0) {11.6}} \put(17,2) {\circle*{0.7}} \put(17,3) {$w$}
\end{picture} \end{center}
    \caption{Second case -- even length.} \label{fig:second}
\end{figure}
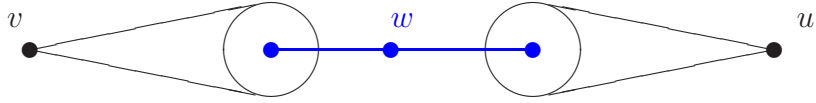
Consider joins of the first type. Recall that we are considering
the process of exposing the neighbourhoods of $u$ and $v$ step by
step. After the first $i-1$ steps, we have exposed the sets
$N_{i-1}(u)$ and $N_{i-1}(v)$, which are assumed to be disjoint.
 Recall that, according to the pairing model (see Introduction), any vertex in $U_{i-1}$ and
$V_{i-1}$ can be regarded as a cell of distinct points, where the
number of points corresponds to the number of unexposed neighbours
of this vertex. The probability that one given point joins another
is then asymptotic to $1/(dn)$, since any pair of unmatched points
is equally likely to be paired and the whole process has, by
Corollary~\ref{wasntused}, matched at most $n^{1/2+o(1)}=o(n)$ pairs
of points to this moment. Asymptotically, there are $|U_{i-1}|
|V_{i-1}|(d-1)^2$ pairs of points such that one is associated with a
vertex in $U_{i-1}$ and the other with a vertex in $V_{i-1}$. This
is because the hypothesis $i = i_0 + k$ implies, by
Lemma~\ref{lem:number_of_nodes}, that the number of vertices in
$U_{i-1}$ or $V_{i-1}$ incident with `bad edges' created at step
$i-2$ is a.a.s.\ at most $O\left(\omega(n) \right)$ for any
$\omega(n) \rightarrow \infty$, and it is clear that each vertex in
$U_{i-1}$ or $V_{i-1}$ with degree larger than 1 in $G[N_{i-1}(u)]
\cup G[N_{i-1}(v)]$ is incident with a `bad edge'.

Thus, the expected number of edges of the first type joining the
neighbourhoods of $u$ and $v$ at time $i-1$, that is, the number of
pairs of points consisting of one point associated with a vertex in
$U_{i-1}$ and one point associated with a vertex of $V_{i-1}$
exposed at time $i$ is asymptotic to $$ \frac{(d-1)^2}{dn} |U_{i-1}|
|V_{i-1}| \,.$$

A similar argument shows that the expected number of edges of the
second type is asymptotic to
$$ \frac{d(d-1)(d-1)^2}{d^2n^2} n |U_{i-1}| |V_{i-1}| = \frac
{(d-1)^3}{dn}|U_{i-1}| |V_{i-1}|\,.$$ Summing these,
\begin{align*}
\E (\textrm{number}& \textrm{ \ of\  joins\  at\  time\  }
i=i_0+k\ |\
S_{i-1}) \\
& \sim \frac{(d-1)^2}{n} |U_{i-1}| |V_{i-1}| \sim
\frac{d^2(d-1)^{2i-2}}{n} = \frac{d^2(d-1)^{2k}}{(d-1)^{2+2
\gamma(n,d)}} .
\end{align*}
We wish to apply the method of moments to establish
$$ \Prob (S_{i_0+k} |\ S_{i_0+k-1})
\sim \exp \left(-d^2(d-1)^{2k- 2\gamma(n,d) -2} \right),$$ so we
have to verify that the $j$-th factorial moment of the random
variable $Z$ counting the number of joins at time $i=i_0+k$
satisfies $\mathbb{E}([Z]_j) = \frac{\mathbb{E}(Z)^j}{j!}$, for
all $j \geq 2$.

Let $j \geq 2$ and suppose that the subgraphs induced by
$N_{i-1}(u)$ and $N_{i-1}(v)$ are disjoint. As before, let $U_{i-1}
= N_{i-1}(u) \setminus N_{i-2}(u)$, $V_{i-1} = N_{i-1}(v) \setminus
N_{i-2}(v)$, and, given $r \in U_{i-1}$, $s \in V_{i-1}$, $t \in
V(G) \setminus (N_{i-1}(v) \cup N_{i-1}(u))$, introduce indicator
random variables $X_{(r,s)}$ for the event that $rs$ is an edge in
$G$ and $Y_{(r,s,t)}$ for the event that $rt$ and $st$ are both
edges in $G$. So,
\begin{equation*}
Z=\sum_{
\begin{array}{lcr}
& \Small{r \in U_{i-1}, ~~ s \in V_{i-1}}
\end{array}
} X_{(r,s)} + \sum_{
\begin{array}{lcr}
& \Small{r \in U_{i-1}, ~~ s \in V_{i-1}}\\
& \Small{t \in V(G) \setminus (N_{i-1}(v) \cup N_{i-1}(u))}
\end{array}
} Y_{(r,s,t)}
\end{equation*}
is the random variable counting the number of joins that appear
between the neighbourhoods of $u$ and $v$ at step $i$.

The $j$-th factorial moment of $Z$ is given by
\begin{equation} \label{eq1}
\mathbb{E}([Z]_j) = \sum_{l=0}^j \sum_\star
\mathbb{P}((X_{(r_m,s_m)}=1,~1 \leq m \leq l) \wedge
(Y_{(r_m,s_m,t_m)}=1,~l+1 \leq m \leq j)),
\end{equation}
where, for any given $l$, $\sum_{\star}$ denotes the sum over all
distinct ordered pairs $(r_m,s_m)$, $1 \leq m \leq l$, and
$(r_{m},s_{m},t_{m})$, $ l+1 \leq m \leq j$.

We shall prove later that the relevant terms in this sum are the
ones for which all the ordered pairs are disjoint, that is, there
is no repetition of vertices among the $j$ events. Assuming this,
we obtain
\begin{equation*}
\begin{split}
\mathbb{E}([Z]_j) = \sum_{l=0}^j \binom{|U_{i-1}|}{j}
\binom{|V_{i-1}|}{j} \binom{n-o(n)}{j-l} & \binom{j}{l}^2 l!
[(j-l)!]^2 \times \\
& \times \Big(\frac{(d-1)^2}{dn-o(n)}\Big)^l
\Big(\frac{(d-1)^3}{dn^2-o(n^2)}\Big)^{j-l}.
\end{split}
\end{equation*}
This is because there are $\binom{|U_{i-1}|}{j} \binom{|V_{i-1}|}{j}
\binom{n-o(n)}{j-l}$ ways of choosing $j$ vertices in each of $U_i$
and $V_i$, and of choosing $j-l$ vertices in $V(G) \setminus
(N_{i-1}(u) \cup N_{i-1}(v))$. Moreover, pairing $l$ of the chosen
vertices in $U_i$ with $l$ of the chosen vertices in $V_i$ can be
done in $\binom{j}{l}^2 l!$ ways, whereas there are $[(j-l)!]^2$
ways of creating triples on the remaining chosen vertices in $U_i$,
$V_i$ and the vertices chosen in $V(G) \setminus (N_{i-1}(u) \cup
N_{i-1}(v))$. Now that we fixed distinct ordered pairs $(r_m,s_m)$,
$1 \leq m \leq l$, and $(r_{m},s_{m},t_{m})$, $ l+1 \leq m \leq j$,
the term $\big(\frac{(d-1)^2}{dn-o(n)}\big)^l
\big(\frac{(d-1)^3}{dn^2-o(n^2)}\big)^{j-l}$ corresponds to the
probability that all the events $X_{(r_m,s_m)}=1$ and
$Y_{(r_m,s_m,t_m)}=1$ occur simultaneously, since there is only a
finite number of them.

The previous sum is asymptotic to
\begin{equation*}
\begin{split}
& \sum_{l=0}^j \frac{|U_{i-1}|^j}{j!} \frac{|V_{i-1}|^j}{j!}
\frac{n^{j-l}}{(j-l)!} \frac{j!^2}{l!}
\Big(\frac{(d-1)^2}{dn}\Big)^l \Big(\frac{(d-1)^3}{dn^2}\Big)^{j-l}\\
& = \frac{|U_{i-1}|^j |V_{i-1}|^j (d-1)^{2j}}{n^j d^j j!}
\sum_{l=0}^j
\frac{(d-1)^{j-l} j!}{l! (j-l)!}\\
& = \frac{|U_{i-1}|^j |V_{i-1}|^j (d-1)^{2j}}{n^j d^j j!}
\sum_{l=0}^j \binom{j}{l} (d-1)^{j-l}\\
& = \frac{|U_{i-1}|^j |V_{i-1}|^j (d-1)^{2j}}{n^j j!} =
\frac{1}{j!} (\mathbb{E}Z)^j =
\frac{1}{j!}\left(\frac{d^2(d-1)^{2k}}{(d-1)^{2+2 \gamma(n,d)}}
\right)^j \;.
\end{split}
\end{equation*}

It remains to show that indeed the sum over all disjoint ordered
pairs $(r_m,s_m)$, $1 \leq m \leq l$, and $(r_{m},s_{m},t_{m})$, $
l+1 \leq m \leq j$, is asymptotic to the sum over all distinct
ordered pairs. Suppose that there are $j-a$ distinct elements
appearing in the first coordinate, $j-b$ in the second and $j-l-c$
in the third, where $a+b+c \geq 1$.  The terms of this form in
equation (\ref{eq1}) are bounded above by
\begin{equation*}
\begin{split}
\sum_{l=0}^j \sum_{**} \binom{j-1}{a} \binom{|U_{i-1}|}{j-a}
&\binom{j-1}{b} \binom{|V_{i-1}|}{j-b}
\binom{j-l-1}{c} \binom{n-o(n)}{j-l-c} \times \\
& \times \binom{j}{l}^2 l! [(j-l)!]^2
\Big(\frac{(d-1)^2}{dn-o(n)}\Big)^l
\Big(\frac{(d-1)^3}{dn^2-o(n^2)}\Big)^{j-l},
\end{split}
\end{equation*}
where $\sum_{**}$ denotes the sum over all triples $(a,b,c) \in
\{0,\ldots,j-1\}^2 \times \{0,\ldots,j-l-1\}$ satisfying $a+b+c \geq 1$.
This is because there are $\binom{|U_{i-1}|}{j-a}$ ways of
choosing $j-a$ vertices in $U_{i-1}$ and $\binom{j-1}{a}$ ways of
building a multi-set of cardinality $j$ with $j-a$  given elements
(and using all of them). The same is true for choosing vertices in
$V_{i-1}$ and $V(G) \setminus (N_{i-1}(u) \cup N_{i-1}(v))$. Our
last expression is smaller or equal to
\begin{equation*}
\begin{split}
\sum_{l=0}^j \sum_{**} \frac{(j-1)^{a+b}(j-l-1)^c}{a!\, b!\, c!} &
\frac{|U_{i-1}|^{j-a}}{(j-a)!} \frac{|V_{i-1}|^{j-b}}{(j-b)!}
\frac{n^{j-l-c}}{(j-l-c)!} \times \\
& \times \frac{j!^2}{l!}\Big(\frac{(d-1)^2}{dn-o(n)}\Big)^l
\Big(\frac{(d-1)^3}{dn^2-o(n^2)}\Big)^{j-l} \;.
\end{split}
\end{equation*}
If we divide this by $ \frac{|U_{i-1}|^j |V_{i-1}|^j
(d-1)^{2j}}{n^j j!}$, this is asymptotic (with respect to $n$) to
\begin{equation*}
\sum_{l=0}^j \sum_{**} \frac{\mathcal{K}(a,b,c,j,l,d)}{|U_{i-1}|^a
|V_{i-1}|^b n^c},
\end{equation*}
where $\mathcal{K}(a,b,c,j,l,d)$ does not depend on $n$. Since
$|U_{i-1}|^a |V_{i-1}|^b n^c \rightarrow \infty$ as $n \rightarrow
\infty$ for every $a+b+c \geq 1$, we conclude that the above sum
tends to zero as $n$ tends to infinity and therefore the terms
related to non-disjoint tuples in equation (\ref{eq1}) can indeed
be ignored to compute $\mathbb{E}([Z]_j)$.

\

Given this, the method of moments implies that
$$ \Prob (\textrm{no\  joins\  at\  time\  } i=i_0+k\  |\  N_{i-1}(v)
\cap N_{i-1}(u) = \emptyset) \sim  \exp
\left(-d^2(d-1)^{2k-2\gamma(n,d)-2}\right),$$ which completes the
proof of the claim and therefore establishes the lemma.
\end{proof}

We are now prepared to prove the result mentioned at the end of
the last section.
\begin{cor}\label{upper bound}
Let $u,v$ be vertices chosen independently with uniform probability
in $G \in \mathcal{G}_{n,d}$. For any function $\omega(n)$ such that
$\omega(n) \to \infty$, the assertion $\dist(u,v) < \log_{d-1}n +
\omega(n)$ holds a.a.s.
\end{cor}

\begin{proof}
Let $\epsilon>0$. Lemma \ref{no join} implies that the probability
of the event $E_k$ that $u$ and $v$ are at distance greater than
$\log_{d-1}n+k$ is asymptotic to $$\exp
\left(-d^2(d-1)^{2k-2\gamma(n,d)-2}\right),$$ where $k$ a fixed
integer and $\gamma(n,d)$ is the fractional part of $\frac{1}{2}
\log_{d-1}n$. So, $\Prob(E_k)<\epsilon$ for $k$ sufficiently large,
and our result follows.
\end{proof}

\begin{lemma}

Let $k$ be an integer and let $\gamma(n,d)$ be the fractional part
of $\frac{1}{2} \log_{d-1} n$. Define $O_i$ to be the random
variable counting the number of $uv$-paths of odd length, that is,
 paths of the first case, created at step $i$. Let
$E_i$ be the equivalent random variable for paths of even length.
Then
\begin{itemize}
\item[(i)]  With $ \mu_k =
d(d-1)^{2k-2\gamma(n,d)-2 }$,
$$\Prob (O_{i_0+k} = j ~|~ N_{i_0+k-1}(u) \cap
N_{i_0+k-1}(v) = \emptyset) \sim \frac{\mu_k^j}{j!}\exp
(-\mu_k).$$

\item[(ii)] With $ {\nu_k =
d(d-1)^{2k-2\gamma(n,d)-1}}$,
$$
\Prob (E_{i_0+k} = j ~|~ N_{i_0+k-1}(u) \cap N_{i_0+k-1}(v) =
\emptyset \wedge O_{i_0+k}=0) \sim \frac{\nu_k^j}{j!}\exp
(-\nu_k).
$$

\end{itemize}

\end{lemma}

\begin{proof}
This can be proven by the method of moments using calculations very
similar to the ones in the previous lemma, proceeding separately for
joins of the first type and joins of the second type. For the
former, we condition on the event that no joins occurred in previous
steps of the process, and, for the latter, we further assume that no
joins of the first type occurred in the current step. The details
are omitted.
\end{proof}
We observe that, alternatively, the proofs of the previous lemma and
of Lemma~\ref{no join} could be unified by considering joint
factorial moments of random variables for joins of the first type
and of the second type.

We are now ready to deduce one of the main results.

\begin{theorem}\label{thm3.5}
Fix an integer $l \geq 1$. The probability that two vertices $u,v$
chosen independently with uniform probability in $G \in
\mathcal{G}_{n,d}$ are joined by exactly $l$ distinct geodesics is
asymptotic to
\begin{equation*}
\begin{split}
\sum_{k=-\infty}^\infty \frac{\big( d(d-1)^{2k-2\gamma(n,d)-2}
\big)^l}{l!} \exp& \left(-\frac{d(d-1)^{2k-2\gamma(n,d)-1}}{d-2}
\right) \times\\
& \times \big(1+(d-1)^l \exp (-d(d-1)^{2k-2\gamma(n,d)-1}) \big)\,.
\end{split}
\end{equation*}
\end{theorem}

\begin{proof}
Let $Z_l$ be the event that $u$ and $v$ are joined by exactly $l$
geodesics, and let $J_i$ denote the event that the first join occurs
at time $i$. Then, given a positive integer $M$,
$$\Prob(Z_l) = \Prob(Z_l \wedge \bigcup_{k=1}^{i_0-M} J_k) + \sum_{k= -
M + 1}^{M-1} \Prob(Z_l \wedge J_{ i_0 +k}) + \Prob(Z_l \wedge
\bigcup_{k  \geq  i_0 + M } J_k).$$ The first and last element in
the right-hand side can be made less than $\epsilon$, for any
given $\epsilon>0$, by choosing $M=M_{\epsilon}$ sufficiently
large, as ensured by Corollary \ref{upper bound} and by the fact
that equation~(\ref{l:first_time}) holds a.a.s. Also, each of the
terms $\Prob(Z_l \wedge J_{ i_0 +k})$, for $-M +1 \leq k \leq
M-1$, is equal to
\begin{equation*}
\begin{split}
\Prob(N_{i_0  +k-1}&(u) \cap N_{i_0  +k-1}(v) = \emptyset) \big[
\Prob(O_{i_0  +k} = l ~|~ N_{i_0 +k-1}(u) \cap N_{i_0 +k-1}(v) =
\emptyset) ~ + \\
& + ~ \Prob(O_{i_0  +k} = 0 ~|~ N_{i_0 +k-1}(u) \cap N_{i_0
+k-1}(v) =
\emptyset) ~ \times \\
& ~ \times \Prob(E_{i_0 +k}=l ~|~ O_{i_0 +k} = 0 \wedge N_{i_0
+k-1}(u) \cap N_{i_0+k-1}(v) = \emptyset) \big].
\end{split}
\end{equation*}
By our previous lemmas, this is asymptotic to
\begin{equation*}
\begin{split}
\exp & \left( -\frac{d(d-1)^{2k-2\gamma(n,d)-2}}{d-2} \right)
\left( \frac{\big( d(d-1)^{2k-2\gamma(n,d)-2}\big)^l}{l!}
\exp (-d(d-1)^{2k-2\gamma(n,d)-2}) ~ \right. + \\
& ~ + \left. \frac{\big( d(d-1)^{2k-2\gamma(n,d)-1}\big)^l}{l!}
\exp (-d(d-1)^{2k-2\gamma(n,d)-2} - d(d-1)^{2k-2\gamma(n,d)-1})
\right)
\end{split}
\end{equation*}

Hence, if we let $\epsilon$ tend to zero,
\begin{equation*}
\begin{split}
\Prob(Z_l) \sim & \sum_{k=-\infty}^\infty  \frac{\big(
d(d-1)^{2k-2\gamma(n,d)-2} \big)^l}{l!} \exp
\left(-\frac{d(d-1)^{2k-2\gamma(n,d)-1}}{d-2} \right)
\times\\
& \times \big(1+(d-1)^l \exp (-d(d-1)^{2k-2\gamma(n,d)-1}) \big),
\end{split}
\end{equation*}
as required.
\end{proof}
An interesting special case is when $l=1$, since this theorem
provides the probability of $u$ and $v$ being joined by a unique
geodesic. This probability is given by
\begin{equation*}
\begin{split}
\sum_{k=-\infty}^\infty d(d-1)^{2k-2\gamma(n,d)-2} \exp &
\left(-\frac{d(d-1)^{2k-2\gamma(n,d)-1}}{d-2} \right)
\times\\
& \times \big(1+(d-1) \exp (-d(d-1)^{2k-2\gamma(n,d)-1}) \big).
\end{split}
\end{equation*}
The probability here is a function of $\gamma(n,d)$ and oscillates
as $\gamma(n,d)$ varies from $0$ to $1$.

We include some numerical results in the table below for some values
of $d$, where $prob$ is the probability of a unique geodesic as
$\gamma(n,d)=0$ and $osc$ is the maximum variation with respect to
$\gamma=0$ as $\gamma$ varies from $0$ to $1$.

\

\begin{center}
\begin{tabular}{|c|c|c|}
\hline $d$ & $prob$ & $osc$\\
\hline $3$ & $0.7213$  & $8.6 \times 10^{-6}$\\
\hline $4$ & $0.6073$  & $1.4 \times 10^{-3}$\\
\hline $5$ & $0.5444$ & $7.9 \times 10^{-3}$\\
\hline $10$ & $0.4411$ & $7.6 \times 10^{-2}$\\
\hline $100$ & $0.3743$ & $0.3$ \\
\hline
\end{tabular}
\end{center}

\

The magnitude of the oscillations depends on $d$. We justify why
it is small when $d$ is small. Note that the probability of a
unique geodesic is equal to
\begin{equation}\label{unique}
\begin{split}
\frac{d-2}{d-1} S_{(d-1)^2}&\left(-\gamma(n,d)+ \log_{(d-1)^2}
\frac{d}{(d-1)(d-2)}\right) + \\
&+ \frac{d-2}{d-1} S_{(d-1)^2}\left(-\gamma(n,d)+ \log_{(d-1)^2}
\frac{d}{d-2}\right),
\end{split}
\end{equation}
where $\displaystyle{S_c(x) = \sum_{m=-\infty}^{\infty}c^{m+x}\exp
(-c^{m+x}) },$ a function with period 1. The classical Poisson
summation formula gives us that
\begin{equation*}
S_c(x) = \sum_{m=-\infty}^{\infty} \int_{-\infty}^{\infty}c^{t+x}
\exp (-c^{t+x}) \exp (2 \pi i m t)\,  dt.
\end{equation*}
Setting $z=c^{t+x}$ gives
\begin{equation}\label{eqPoisson} S_c(x)
= \frac{1}{\log c}\sum_{m=-\infty}^{\infty} \exp(-2 \pi i mx)
\int_{0}^{\infty} \exp(-z + 2 \pi i m \log z/ \log c)\, dz,
\end{equation}
and the integral is just $\Gamma(2 \pi i m/ \log c + 1)$.

By properties of the gamma function (see for instance
\cite{numerics}), we have
\begin{equation*}
|\Gamma(1 + yi )| = |iy \Gamma(yi)| = |y| \sqrt{\frac{\pi}{y
 \sinh(\pi y)}},
\end{equation*}
so given $m$ in the previous summation,
\begin{equation}\label{bound}
\begin{split}
|\Gamma(2 \pi m i/ \log c + 1)\exp(-2 \pi i mx)| &= |\Gamma(2 \pi
i m/ \log c + 1)|\\
&= \left( \frac{2 \pi^2|m|}{\log{c} ~ |\sinh(2 \pi^2m/\log{c})|}
\right)^{1/2}.
\end{split}
\end{equation}

The term for $m=0$ in the sum (\ref{eqPoisson}) is independent of
$x$, hence it yields terms independent of $n$ in equation
(\ref{unique}). In the special case $d=3$, equation (\ref{bound})
leads to the following bounds on the other terms of the sum
(\ref{eqPoisson}). For $|m|=1$, the bound is approximately $4.32
\times 10^{-3}$, for $|m|=2$, it is approximately $4.94 \times
10^{-6}$, and for larger $|m|$ the bounds are even smaller.
Similar observations explain the small oscillations when $d$ is
small.

\

\section{Almost geodesic cycles}\label{almost}

In the proof of Theorem~\ref{main} we shall use the following
auxiliary result.

\begin{lemma}\label{l:aux}
Let $G \in \mathcal{G}_{n,d}$ and let $u,v$ be vertices chosen
independently at random in $G$. Consider functions
$\alpha(n),\beta(n)$ such that $\alpha(n),\beta(n) \rightarrow
\infty$, $\alpha(n) = o(\log_{d-1}n)$  and $\beta(n) =
o(\alpha(n))$. Then a.a.s.\ every vertex at distance $\lfloor
\alpha(n) \rfloor$ from $u$ or $v$ lies on at most one $uv$-path
with length less than or equal to $\log_{d-1}n + \beta(n)$.
\end{lemma}

\begin{proof}
We prove this result for vertices at distance $\lfloor \alpha(n)
\rfloor$ from $u$, and by a similar argument the same result holds
for vertices at distance $\lfloor \alpha(n) \rfloor$ from $v$. As
in Section \ref{section3}, we consider the process of exposing the
neighbourhoods of $u$ and $v$   based on the pairing model. Here,
$N(u)$, the neighbourhood of $u$, is exposed for $\lfloor
\alpha(n) \rfloor$ steps while $N(v)$, the neighbourhood of $v$,
is exposed for $\lfloor \frac{1}{2} \log_{d-1}n  - \beta(n)
\rfloor$ steps. By Lemma \ref{lem:number_of_nodes}, a.a.s.\ $N(u)$
and $N(v)$ are both trees. Moreover, Lemma \ref{no join} ensures
that $N(u) \cap N(v) = \emptyset$ holds a.a.s.

Let $U_{\alpha}$ be the set of vertices at distance $\lfloor
\alpha(n) \rfloor$ from $u$. Given a vertex $w \in U_\alpha$, let
$Y_w$ be the indicator random variable for the event that $w$ is
connected to $N(v)$ by at least two distinct paths of length less
than or equal to $ \frac{1}{2} \log_{d-1}n  - \alpha(n) + 2\beta(n)
+ 2$. Define $Y = \sum_{w \in U_{\alpha}} Y_i$. It is clear that
this lemma follows if we prove that a.a.s.\ $Y=0$. We shall do this
by using $$\Prob \big( Y \geq 1 \big) \leq \sum_{w \in U_{\alpha}}
\Prob \big( Y_w = 1 \big),$$ and by showing that the right-hand side
goes to zero as $n$ tends to infinity.

For a fixed $w$, define the set $N_w^{\prime}$ obtained by the
exposure of the neighbourhood of $w$ for $ \frac{1}{2} \log_{d-1}n -
\alpha(n) + 2\beta(n) + 2$ steps. This time, however, the neighbour
of $w$ in $N(u)$ is not added to $N_w^{\prime}$ at the first step of
the process, that is, only the ``new" neighbours of $w$ are exposed.
As in Lemma \ref{lem:number_of_nodes}, we use the term ``bad edges"
for edges that yield cycles in $N_w^{\prime}$. Consider the random
variable $X_w$ counting the number of ``bad edges" in
$N_w^{\prime}$. Then, calculations analogous to the ones in Lemma
\ref{lem:number_of_nodes} establish that
\begin{equation*}
\begin{split}
\mathbb{E}(X_w) &= \sum_{s=0}^{\lfloor \frac{1}{2} \log_{d-1}n  -
\alpha(n) + 2\beta(n) + 2 \rfloor} O\left( \frac{(d-1)^{2s}}{n} \right) \\
& = O \left( \frac{(d-1)^{\log_{d-1}n  - 2\alpha(n) +
4\beta(n)}}{n} \right) = O\left( (d-1)^{4\beta(n)- 2\alpha(n)}
\right).
\end{split}
\end{equation*}
Thus Markov's inequality implies $$\Prob\big(X_w \geq 1\big) =
O\left( (d-1)^{4\beta(n)- 2\alpha(n)} \right).$$ Now, note that
$$\Prob \big(Y_w = 1 \big) = \Prob \big(X_w \geq 1 \big) \Prob
\big(Y_w = 1 ~|~ X_w \geq 1\big) + \Prob \big(X_w = 0 \big) \Prob
\big(Y_w = 1 ~|~ X_w = 0 \big).$$ We have a bound for the first term
in this sum. For the second term, we use the definition of
conditional probability and observe that the event $(Y_w = 1) \wedge
(X_w=0)$ occurs only if there is a pair of distinct paths joining
$w$ to $N(v)$ with length at most $\frac{1}{2} \log_{d-1}n -
\alpha(n) + 2\beta(n) + 2$ and with the property that, after they
first split, they do not join again.

So, a bound on $\Prob \big(Y_w = 1 \wedge X_w = 0 \big)$ may be
obtained by counting the number of possible pairs of distinct paths
$P$ and $Q$ joining $u_i$ to $N(v)$ with lengths $r$ and $s$, where
$r \leq s \leq \frac{1}{2} \log_{d-1}n - \alpha(n) + 2\beta(n) + 2$,
and the first $j$ vertices are shared by both paths, while the
remainder of the paths are internally disjoint. So, if $i_0 =
\lfloor \frac{1}{2} \log_{d-1} n \rfloor$,
\begin{equation*}
\begin{split}
\Prob \big(Y_w = 1 \wedge X_w = 0 \big) &= \sum_{s=1}^{\lfloor i_0
- \alpha(n) + 2\beta(n) + 2\rfloor} \sum_{r=1}^{\lfloor i_0 - s
\rfloor} \sum_{j=0}^{r-1}\\
&~~~~\binom{(d-1)^{\lfloor i_0-\beta(n) \rfloor}}{2}
\binom{n-o(n)}{r+s-j-2}
\binom{r+s-j-2}{j}\binom{r+s-2}{r-1}\\
&~~~ \times j!
(r-1)! (s-1)! O\left(\left(\frac{(d-1)}{n-o(n)}\right)^{r+s-j} \right)\\
& = \sum_{s=1}^{\lfloor i_0 - \alpha(n) + 2\beta(n) +2 \rfloor}
\sum_{s=1}^{\lfloor i_0 - s \rfloor} \sum_{j=0}^{r-1} O\left( \frac{(d-1)^{2i_0-2 \beta(n)}(d-1)^{r+s-j}}{n^2} \right)\\
&= O \left((d-1)^{2\beta(n)- 2\alpha(n)} \right).
\end{split}
\end{equation*}
Note that the formula holds because there are at most
$\binom{(d-1)^{\lfloor i_0-\beta(n) \rfloor}}{2}$ ways of choosing
two vertices in $N(v)$ and there are $\binom{n-o(n)}{r+s-j-2}$ ways
of choosing vertices in the graph to include in the two paths.
Moreover, these vertices can be divided into vertices of $P \cap Q$,
$P-Q$ and vertices of $Q-P$ in
$\binom{r+s-j-2}{j}\binom{r+s-2}{r-1}$ ways and can then be ordered
to form the paths in $j! (r-1)! (s-1)!$ ways. Finally, each edge on
the path appears with probability at most $\frac{(d-1)}{n-o(n)}$
conditional on the fact that all previous edges on the path have
appeared.

We conclude that
\begin{equation*}
\begin{split}
\Prob \big(Y_w = 1 \big) &= \Prob \big(X_w \geq 1 \big) \Prob
\big(Y_w = 1 ~|~ X_w \geq 1\big) + \Prob \big(X_w = 0 \big) \Prob
\big(Y_w = 1 ~|~ X_w = 0 \big)\\
& \leq \Prob \big(X_w \geq 1 \big) + \Prob \big(Y_w = 1 ~|~ X_w =
0 \big)\\
&= O \left((d-1)^{4\beta(n)- 2\alpha(n)} \right) +
\frac{O\left((d-1)^{2\beta(n)- 2\alpha(n)}
\right)}{1-O\left((d-1)^{4\beta(n)- 2\alpha(n)} \right)}\\
& = O \left((d-1)^{4\beta(n)- 2\alpha(n)} \right).
\end{split}
\end{equation*}
Now, because there are $O\left( (d-1)^{\alpha(n)} \right)$
vertices at distance $\lfloor \alpha(n) \rfloor$ of $u$, we have
\begin{equation*}
\begin{split}
\Prob \big(Y \geq 1 \big) & \leq \sum_{w \in U_{\alpha}} \Prob
\big( Y_w = 1
\big)\\
& = O\left( (d-1)^{\alpha(n)} \right)O \left((d-1)^{4\beta(n)-
2\alpha(n)} \right) = O \left((d-1)^{4\beta(n)- \alpha(n)} \right).
\end{split}
\end{equation*}
Because $\beta(n)=o(\alpha(n))$, this term goes to zero as $n$
tends to infinity and indeed $\Prob \big(Y \geq 1 \big)
\rightarrow 0$. The lemma follows.
\end{proof}

We are now ready to prove the main theorem.
\begin{proof}[Proof of Theorem~\ref{main}]
Part (i) of the theorem follows from Corollary \ref{upper bound}
and equation (\ref{l:first_time}).

Before proving part (ii), it is convenient to deal first with the
following simpler goal.  To state this we need two definitions. A
{\em $k$-near-geodesic} is a path  that is a geodesic between the
two vertices at distance $k$ from its ends. A vertex $p$ on a path
$P$ between vertices $u$ and $v$ is said to be a {\em midpoint} of
$P$ if $|\dist_{P}(u,p)-\dist_{P}(v,p)| \leq 1$, where $\dist_{P}$
denotes the distance along path $P$.

\begin{lemma}\label{waspartii}
Asymptotically almost surely, for every  two distinct
$k$-near-geodesics, $P$ and $Q$, between $u$ and $v$ with
midpoints $p$ and $q$, respectively,
$$
\log_{d-1}(n) - \omega(n) < \dist(p,q) < \log_{d-1}(n) +
\omega(n).
$$
\end{lemma}
\begin{proof}
For the upper bound, we know as in part (i) of the theorem  that
a.a.s.\ $\dist(u,v) < \log_{d-1}n + \omega(n)$, hence there is a
sufficiently short path connecting $p$ to $q$ through $u$ or $v$.

We turn to the lower bound. Given a function $\omega(n)$ satisfying
$\omega(n) \rightarrow \infty$, we know that a.a.s.\ $\dist(u,v)
\geq \log_{d-1}n -  \omega(n)$ (see~(\ref{l:first_time})). Consider
distinct $k$-near geodesics $P$ and $Q$.

\noindent {\bf Claim 1:\ } {\em $P$ and $Q$ a.a.s.\ do not have a
vertex in common at distance at least $\frac{\omega(n)}{3}$ from
their endpoints.}

We prove the claim by contradiction. Suppose without loss of
generality that such a vertex is closer to $u$ on $P$ and let $w$ be
the vertex on $P$ at distance $\lfloor \frac{\omega(n)}{3} \rfloor$
from $u$. Note that $P$ and $Q$ a.a.s.\  differ at some vertex or
edge after $w$, since the set of vertices at distance at most
$\frac{\omega(n)}{3}$ a.a.s.\ induces a tree. But then, $w$ lies on
at least two distinct $u,v$ paths with length less than or equal to
$\log_{d-1}n + \log \omega(n)$, which a.a.s.\ does not occur by
Lemma \ref{l:aux} with $\alpha(n)=\frac{\omega(n)}{3}$ and
$\beta(n)=\log \omega(n)$. (Note that the lengths of both $P$ and
$Q$ are a.a.s.\ bounded by $\log_{d-1}n + \log \omega(n)$ because
any $k$-near geodesic between $u$ and $v$ has length at most
$\dist(u,v)+4k$.) This proves the claim.

Now consider the event that the midpoints $p$ and $q$ of  $P$ and
$Q$ are at distance at most $\log_{d-1}n - \omega(n)$. One way for
this to occur is by the existence of a $pq$-path $R$ of length at
most $\log_{d-1}n - \omega(n)$ using vertices and edges on $P \cup
Q$ only.  But $\dist(u,v) \geq \log_{d-1}n -  \omega(n)$  implies
that $R$ does not contain vertices at distance less than or equal to
$\frac{\omega(n)}{3}$ from $u$ or $v$.   Claim~1   shows that no
other vertex can be in common. Thus, a.a.s.\ there is no short path
from $p$ to $q$ using edges on $P$ and $Q$ only.

So consider a geodesic $A$ between $p$ and $q$ containing at least
one edge outside $P \cup Q$. Using $A$ oriented from $p$ to $q$ as a
reference, let $v_P$ denote the last vertex on $A \cap P$ and let
$v_Q$ be the first vertex on $A \cap Q$ after $v_P$. The vertices
$v_P$ and $v_Q$ divide the geodesic into three parts, namely from
$p$ to $v_P$, from $v_P$ to $v_Q$ and from $v_Q$ to $q$. Because
$P,Q$ are $k$-near-geodesics between $u$ and $v$ for a fixed $k$ and
$A$ is a geodesic between $p$ and $q$, we must have $\dist_A(p,
v_P)=\dist_P(p,v_P)$ and $\dist_A(v_{Q},q)=\dist_{Q}(v_{Q},q)$. So,
for $p$ and $q$ to be at distance at most $\log_{d-1}n - \omega(n)$
for some $\omega(n) \rightarrow \infty$, it must be that
\begin{equation}\label{Aeq}
\dist_A(v_P,v_Q) < \log_{d-1}n - \dist_P(p, v_P) -
\dist_{Q}(v_Q,q) - \omega(n).
\end{equation}
So, a short path between $p$ and $q$ has to be caused by a short
path connecting a vertex in $P$ to a vertex in $Q$ which is
internally disjoint from $P \cup Q$. More precisely, there must
exist vertices $v_P,v_Q$ on $P$ and $Q$, respectively, and an
$v_Pv_Q$-path $R$ satisfying:
\begin{eqnarray}\label{eq2}
V(R) \cap V(P \cup Q) = \{v_P,v_Q\},\\
|R| \leq \log_{d-1}n - \dist(p,v_P) - \dist(v_Q,q) -
\omega(n).\label{other}
\end{eqnarray}
\begin{figure}
\begin{center}
\setlength{\unitlength}{0.5cm}
\begin{picture}(26,4)
\color{red} \thicklines \put(0,2) {\line(2,1){4}} \put(4,4)
{\line(2,0){2}} \put(8,4){\line(2,0){4}} \put(14,4){\line(2,0){4}}
\put(20,4){\line(2,0){2}} \put(22,4){\line(2,-1){4}}
\put(24,3.5){$P$} \multiput(6,4)(0.2,0){10} {\line(1,0){0.1}}
\multiput(12,4)(0.2,0){10} {\line(1,0){0.1}}
\multiput(18,4)(0.2,0){10} {\line(1,0){0.1}} \color{blue}
\thicklines \put(0,2) {\line(2,-1){4}} \put(4,0) {\line(2,0){2}}
\put(8,0) {\line(2,0){4}} \put(14,0){\line(2,0){4}}
\put(20,0){\line(2,0){2}} \put(22,0) {\line(2,1){4}}
\put(24,0){$Q$} \multiput(6,0)(0.2,0){10}{\line(1,0){0.1}}
\multiput(12,0)(0.2,0){10}{\line(1,0){0.1}}
\multiput(18,0)(0.2,0){10}{\line(1,0){0.1}} \color{black}
\multiput(4,0)(2,0){10} {\circle*{0.2}} \multiput(2,1)(22,0){2}
{\circle*{0.2}} \multiput(0,2)(26,0){2} {\circle*{0.2}}
\multiput(2,3)(22,0){2} {\circle*{0.2}} \multiput(4,4)(2,0){10}
{\circle*{0.2}} \put(0,2.5) {$u$} \put(26,2.5) {$v$}
\put(10,4.5){$v_P$} \put(16,4.5){$p$} \put(10,0.5){$q$}
\put(16,0.5){$v_Q$} \color{green} \put(10,4){\line(3,-2){2}}
\put(16,0){\line(-3,2){2}}
\multiput(12,2.66)(0.2,-0.133){10}{\line(1,0){0.1}}
\put(14,2){$R$}
\end{picture}
\end{center}
      \label{fig:fourth_form}
      \caption{Path $R$}
\end{figure}
Such a configuration is illustrated by Figure 3.

We prove that a.a.s.\ $G$ does not contain a path $R$ satisfying
(\ref{eq2}) and (\ref{other}). We do this by exposing the
neighbourhoods of vertices along $P$ and $Q$ conditional on the
particular paths $P$ and $Q$ being in the graph. By relaxing the
condition that $P$ and $Q$ are $k$-near geodesics (but retaining the
condition that their length is at most $\dist(u,v)+4k+\omega(n)$),
we may take the rest of the pairing to be random. We will later
argue that the number of pairs of such paths $P$ and $Q$ is small
enough for our argument to work.

Given a vertex $p_r$ at distance $r$ from $p$ along $P$ and a vertex
$q_s$ at distance $s$ from $q$ along $Q$, let $X_{p_r,q_s}$ be the
event that $p_r$ and $q_s$ are connected by a path of length at most
$\log_{d-1}n - r - s - \omega(n)$ which is internally disjoint from
$P$ and $Q$. Define the random variable $Y_{P,Q}=\sum_{r=0}^{\lceil
\log_{d-1}n - \omega \rceil} \sum_{s=0}^{\lceil \log_{d-1}n - \omega
\rceil - r} X_{p_r,q_s}$, so that $Y_{P,Q}=0$ only if $G$ does not
contain a path $R$ satisfying (\ref{eq2}) and (\ref{other}) with
respect to $P$ and $Q$.

Once again, we look at the process in which the neighbours of $p_r$
and $q_s$ are exposed, then their neighbours are exposed, and so on,
but we do not consider the neighbours of $p_r$ and $q_s$ on $P$ or
$Q$, so as to expose the sets $N_i(p_r)$ and $N_i(q_s)$ containing
only the vertices at distance at most $i$ from $p_r$ and $q_s$ that
can be reached by paths internally disjoint from $P$ and $Q$.
Clearly, $p_r$ and $q_s$ are joined by a path as in (\ref{eq2}) and
(\ref{other}) only if $N_i(p_r)$ and $N_i(q_s)$ join in at most
$\frac{1}{2}(\log_{d-1}n - r - s - \omega(n))$ steps. The
probability of this can be calculated as in the earlier sections,
and we conclude that
\begin{equation*}
\begin{split}
\Prob(Y_{P,Q} \geq 1) & \leq \E(Y) \\
& \leq \sum_{r=0}^{\lceil \frac{1}{2}(\log_{d-1}n - \omega)
\rceil}
\sum_{s=0}^{\lceil \frac{1}{2}(\log_{d-1}n - \omega \rceil) - r} \Prob(X_{p_r,q_s})\\
& = 4 \sum_{r=0}^{\lceil \frac{1}{2}( \log_{d-1}n - \omega)
\rceil}
\sum_{s=0}^{\lceil \frac{1}{2}(\log_{d-1}n - \omega )\rceil- r} O\left((d-1)^{2\left(\frac{1}{2}(\log_{d-1}n - \omega) -r -s\right)}/n\right)\\
& = O\left((d-1)^{-\omega(n)}\right).
\end{split}
\end{equation*}
By Lemma \ref{l:aux} with $\alpha(n)$ any function tending to
infinity sufficiently slowly, and $\beta(n)=o(\alpha(n))$, the
number paths of length at most $\log_{d-1}n + \beta(n)$ between $u$
and $v$ is a.a.s.\ at most $ 2(d-1)^{\alpha(n)}$, since this is a
bound on the number of vertices at distance $\lfloor \alpha(n)
\rfloor$ from $u$. Thus, a.a.s.\ there are at most $ \gamma(n) $
pairs of such paths between $u$ and $v$, for any
$\gamma(n)\to\infty$. Let $Z$ denote this asymptotically almost sure
event (for some $\alpha$ to be restricted shortly), and let
$\bar{Z}$ be its complement.

Let $Y=\sum_{P,Q} Y_{P,Q}$, where the sum is over all pairs of
distinct paths between $u$ and $v$ whose length is at most
$\dist(u,v)+4k+\omega(n)$. Then,
\begin{equation*}
\begin{split}
\Prob(Y \geq 1) &= \Prob((Y \geq 1) \cap \bar{Z}) + \Prob((Y \geq
1) \cap
Z)\\
& \leq  \Prob(\bar{Z}) + \sum_{P,Q} \Prob((Y_{P,Q} \geq 1)\cap
Z)\\
& \leq \Prob(\bar{Z}) + O\left(\gamma(n)(d-1)^{-\omega(n)}\right),
\end{split}
\end{equation*}
which tends to $0$ provided $\gamma(n)= o((d-1)^{\omega(n)})$. As a
consequence, a.a.s.\ there are no configurations satisfying the
conditions in (\ref{eq2}) and (\ref{other}). Hence, a.a.s.\ the
inequality $\dist(p,q) > \log_{d-1}n - \omega(n)$ holds, and we the
lemma.
\end{proof}

We return to the proof of the theorem. We say that two points on a
cycle are {\em diametrically opposite} if the distance between them
around the cycle is $\lfloor \ell/2\rfloor$, where the cycle has
length $\ell$. Note that  part (ii) of the theorem follows
immediately from the following, since if there is a ``short-cut" for
any two vertices on a cycle, there is a short-cut for a pair of
diametrically opposite ones.

\noindent {\bf Claim 2:\ } {\em  Asymptotically almost surely, there
is a cycle $C$  passing through $u$ and $v$, of length $\ell$
satisfying $|2\log_{d-1}n - \ell|< \omega(n)$, such that for every
pair of points $p$ and $q$ that are diametrically opposite on $C$,
$d(p,q)\ge \log_{d-1}n-\log_{d-1}\log_{d-1}n-\omega(n)$.}

We now prove Claim~2. Note that, by (i), for all $\eps>0$, there is
$K$ sufficiently large that $\pr(d(u,v)>\log_{d-1}n -K)>1-\eps$ for
all $n$ sufficiently large. From the proof of Lemma~\ref{no join} it
is evident that for each neighbour $u'$ of $u$ and neighbour $v'$ of
$v$, there is a.a.s.\ a path from $u$ to $v$ in which the second and
penultimate vertices are $u'$ and $v'$, and of length at most
$\log_{d-1}n +\omega(n)$.  Hence we may also assume that $K$ is
sufficiently large that the probability that there is such a path
for two given neighbours $u'$ and $v'$ and of length at most
$\log_{d-1}n +K$ is at least $1-\eps$.  So with probability at least
$1-2\eps$, we can choose two such $uv$-paths $P$ and $Q$, where the
neighbours of $u$ and $v$ on $P$ are both different from those on
$Q$. In each case we may select a shortest path with these
specifications. Then $P$ and $Q$ must be $(K+1)$-near geodesics.

By Claim~1, we may assume there is no vertex in common between $P$
and $Q$ that is more than $\omega(n)/3$ from $u$ and $v$. For
$\omega$ growing slowly enough, there is a.a.s.\ no point in common
that is at most $\omega(n)/3$ from $u$ and $v$ either, since
Lemma~\ref{lem:number_of_nodes} implies that a.a.s.\ neither $u$ nor
$v$ is in a short cycle. Let $C$ be the union of the paths $P$ and
$Q$. From the bounds on $d(u,v)$, $C$ has length at least
$2\log_{d-1}n -2K$.

To prove the statement about all diametrically opposite points $p$
and $q$ on $C$, we may rework the argument in Lemma~\ref{waspartii}.
The Claim proved above shows that every short path of the type we
are interested in must use some edge not on $P$ or $Q$. Arguing as
before, we only need to eliminate the existence of $A$ such
that~(\ref{Aeq}) holds. The same argument as before shows that for
any fixed such $p$ and $q$, with $Y_{P,Q}$ defined as before, we
again have $\Prob(Y_{P,Q} \geq 1) = O\left((d-1)^{-
\omega(n)}\right)$.

Now apply this inequality to the $O(\log_{d-1}n)$ pairs of vertices
$p$ and $q$ diametrically opposite on $C$. Also, put
$$
\omega(n)= \log_{d-1}\log_{d-1} n +\gamma(n).
$$
Then the probability that $Y_{P,Q} \geq 1$ for at least one of these
choices of $p$ and $q$ is $O\left((d-1)^{- \gamma(n)}\right)$.
Hence, if $\gamma(n)\to\infty$,  we have a.a.s.\ for all such $p$
and $q$,   $d(p,q)\ge f(n)-\log_{d-1}\log_{d-1} n -\gamma(n)$.
Replacing $\gamma$ by $\omega$ gives the final statement in Claim~2,
with probability at least $1-2\eps+o(1)$. This statement is true for
all $\eps>0$. That fact implies that the final statement in Claim~2
holds a.a.s. (This can be regarded as ``letting $\eps\to 0$
sufficiently slowly''.) Combining this with part (i) proves Claim~2,
since, although there may be different functions at the different
occurrences of $\omega$, they can be made the same. This completes
the proof of Claim~2. To obtain the theorem, we note that the proof
of Claim~2 does more: it shows that the two paths can be chosen to
use distinct neighbours of $u$ in their initial step, and distinct
neighbours of $v$ in their final step.
\end{proof}

\section{Final remarks} \label{final}
In this article we have examined the ``shape" of random regular graphs. This brings up related questions.

Our proof of the main theorem can be seen to give more: a.a.s.\ for
every pair of short (i.e.\ bounded length) paths, one containing
$u$ and one containing $v$, there is an almost geodesic cycle
containing both of these paths. We also show that the paths referred
to in the theorem each contain a geodesic between the two vertices
at distance $k$ from its ends, for any $k$ tending to infinity  with
$n$.

A geodesic cycle $C$ in $G$ is a cycle in which for every two
vertices $u$ and $v$ in $C$, the distance $d_G(u,v)$ is equal to
$d_C(u,v)$. A significant open problem is to determine whether in a
random $d$-regular graph, a.a.s.\ almost all pairs of vertices lie
in a geodesic cycle. It is not even known if at least one geodesic
cycle of length asymptotic to $\log_{d-1}n$ exists a.a.s.

We may also draw conclusions on how ``thin" the topological
triangles are in random regular graphs. Consider  the proof of
Lemma~\ref{no join}, which analyses the time at which two
simultaneous breadth-first reaches from $u$ and from $v$ join each
other. The proof is concerned with an accurate estimate of the
probability that there are no joins by a time near $i_0$. It is
easy to see from the ideas in the proof that for large $K$, the
second join is quite likely to occur by time $i_0+K$, and
furthermore that the first two joins are quite likely to be in
branches that diverged, in the breadth first search from $u$,  at
time less than $K$, and similarly from  $v$. Let $u'$ and $v'$ be
the points of divergence near $u$ and $v$. Then the joins give two
paths $P$ and $Q$ from $u'$ to $v'$, the shorter of which, say $P$, is
geodesic, and   we can choose another vertex, $w$, on $Q$, of
distance $K$ from $u'$, such that the resulting two subpaths of
$Q$ to $u'$ and $v'$ from $w$ are both geodesic. Thus $u',v',w$
form a geodesic triangle. By Lemma~\ref{waspartii} (noting $P$ and
$Q$ are $2K$-near geodesics from $u$ to $v$), the distance between
the midpoints of $P$ and $Q$ is a.a.s.\ at least $\log_{d-1}(n) -
\omega(n)$, where $\omega(n)$ is any function tending to $\infty$.
Hence the midpoint $p$ of $P$ a.a.s.\ has distance at least
$\frac12\log_{d-1}(n) - \omega(n)$ from the union of the other two
sides of the geodesic triangle $u',v',w$. The probability in the
above statements tends to 1 if we let $K\to\infty$ sufficiently
slowly. It follows that a random $d$-regular graph, for $d\ge 3$,
is a.a.s.\ $\delta$-hyperbolic for $\delta = (\log_{d-1} n)/2 -
\omega(n)$. An upper bound on such $\delta$ is obviously half of
the diameter of the graph, which is   $(\log_{d-1} n)/2 + O(\log
\log n)$ by the main result of~\cite{BF}.

Finally, it would be interesting to see to what extent the
geometric properties we have addressed in this paper are preserved
if the model of regular graphs changes. Particular models of
interest might be random Cayley graphs, random lifts of regular
graphs, and one can consider also some deterministic models of
expanders.



\end{document}